\input amstex
\documentstyle{amsppt}
\input bull-ppt
\keyedby{bull286e/PAZ}

\define\bbc{{\Bbb C}}
\define\bbd{{\Bbb D}}

\define\bbp{{\Bbb P}}

\define\bbz{{\Bbb Z}}

\define\grO{\Omega}

\define\gri{\iota}

\define\gro{\omega}

\define\la#1{\hbox to #1pc{\leftarrowfill}}

\define\fract#1#2{\raise4pt\hbox{$ #1 \atop #2 $}}
\define\decdnar#1{\phantom{\hbox{$\scriptstyle{#1}$}}
\left\downarrow\vbox{\vskip15pt\hbox{$\scriptstyle{#1}$}}%
\right.}
\define\decupar#1{\phantom{\hbox{$\scriptstyle{#1}$}}
\left\uparrow\vbox{\vskip15pt\hbox{$\scriptstyle{#1}$}}%
\right.}
\define\hsmash{\triangleright\kern-4.4pt\raise.66pt\hbox{$%
\scriptstyle{<}$}}
\define\boxit#1{\vbox{\hrule\hbox{\vrule\kern3pt
\vbox{\kern3pt#1\kern3pt}\kern3pt\vrule}\hrule}}

\topmatter
\cvol{26}
\cvolyear{1992}
\cmonth{April}
\cyear{1992}
\cvolno{2}
\cpgs{317-321}
\title The Atiyah-Jones Conjecture\endtitle
\author C. P. Boyer, J. C. Hurtubise, 
B. M. Mann, and R. J. Milgram\endauthor
\address Department of Mathematics and Statistics,
University of New Mexico, Albuquerque, New Mexico 
87131\endaddress
\ml cboyer\@gauss.unm.edu\endml
\address Department of Mathematics and Statistics,
McGill University, Montr\'eal, Qu\'ebec H3G 1M8, 
Canada\endaddress
\ml hurtubis\@gauss.math.mcgill.ca\endml
\address Department of Mathematics and Statistics,
University of New Mexico, Albuquerque, New Mexico 
87131\endaddress
\ml  mann\@gauss.unm.edu\endml
\address Department of Mathematics,
Stanford University, Stanford, California 94305\endaddress
\ml milgram\@gauss.stanford.edu\endml
\thanks During the preparation of this work
the first and third authors were supported 
by NSF grants, the second author by an NSERC grant, and
the fourth author by MSRI and an NSF grant\endthanks
\date August 14, 1991\enddate
\abstract The purpose of this note is to announce our proof
of the Atiyah-Jones conjecture concerning the topology of 
the
moduli spaces of based ${\scriptstyle SU(2)}$-instantons 
over 
${\scriptstyle S^4}$.   Full details and proofs appear in 
our paper
[BHMM1].\endabstract 
\endtopmatter

\document

\heading 1. Introduction\endheading

In 1978 Atiyah and Jones [AJ] proved a 
foundational theorem on the 
low-dimensional homology of the moduli space of based 
$SU(2)$
instantons on $S^4$ and made a number of conjectures.  In 
this note
we briefly explain how, using methods of algebraic 
geometry and
algebraic topology, we analyze the global geometry of 
these moduli
spaces, which we denote by ${\cal M}_k$, and prove the 
Atiyah-Jones
conjectures.  More precisely, let $P_k$ be a principal
$SU(2)$ bundle over $S^4$ with second Chern class $k$. Let 
${\cal M}_k$ 
be the space of  based
gauge equivalence classes of  connections on $P_k$ that 
satisfy
the Yang-Mills self-duality equations of $SU(2)$ gauge 
theory. 
There is a natural forgetful map $\vartheta_k: {\cal M}_k
\to {\cal B}_k$ where ${\cal B}_k$ is the space of based
equivalence classes of all connections in $P_k$.  Atiyah 
and Jones
[AJ] showed that ${\cal B}_k$ is homotopy equivalent to 
$\grO^3_k
SU(2)$ and proved

\thm{Theorem \rm [Atiyah-Jones]}
$\vartheta_k$
induces a homology surjection 
$$ H_t({\cal M}_k) \buildrel{(\vartheta_k)_t}
\over{\to}
H_t(\grO^3_k SU(2)) \to 0$$
for $t\leq q=q(k)\ll k$. 
\ethm

The main point of the Atiyah-Jones result is that the 
homology of
$\grO^3 SU(2)$ is completely known so that the surjection 
sheds
some light, for small values of $t$ relative to $k$, on the
heretofore unknown groups $H_t({\cal M}_k)$.  Notice that, 
as $k$
increases, the target spaces $\grO^3_k SU(2) \cong
\grO^3_{k+1}SU(2)$ remain homotopy equivalent, whereas the 
topology
of the ${\cal M}_k$, which are $4k$ dimensional complex 
manifolds,
changes as $k$ varies.  Atiyah and Jones then made the 
following
conjectures: 
\roster
\item"1."  $(\vartheta_k)_t$ is a homology isomorphism for 
$t\leq
q(k)$.
\item"2."  The range of the surjection (isomorphism) $q = 
q(k)$ 
can be explicitly determined as a function of $k$.
\item"3."  The homology statements can be replaced by 
homotopy
statements in both conjectures 1 and 2. 
\endroster

The last and strongest statement became commonly known as 
the
Atiyah-Jones conjecture. Notice that, while it is a simple 
fact
that there are natural maps $j_k:\grO^3_k SU(2)
\buildrel{\simeq}\over{\rightarrow}
 \grO^3_{k+1}SU(2),$
there is no obvious map on the ${\cal M}_k$ level that 
raises the
instanton number $k$. Rather, it follows from analytic 
results
of Taubes [T1] that there are compatible inclusions
$\gri_k : {\cal M}_k \to {\cal M}_{k+1}$
on the instanton moduli space level such that the 
following diagram
homotopy commutes
$$\matrix
{\cal M}_k&\buildrel{\gri_k}\over{\rightarrow}&{\cal M}_{k+
1}\\
\decdnar{\vartheta_k}&&\decdnar{\vartheta_{k+1}}\\
\grO^3_k SU(2)& \buildrel{\simeq}\over{\rightarrow}& 
\grO^3_{k+1}SU(2)
\endmatrix
\tag1.1$$
This diagram permits one to take direct limits and hence 
formulate
the stable version of the Atiyah-Jones conjecture.  This 
was proven
by Taubes [T2].

\thm{Theorem \rm [Taubes Stability Theorem]}
Let ${\cal M}_{\infty}$ be the direct limit of the 
${\cal M}_k$\RM's under the inclusions $\gri_k$ and let
$\vartheta:{\cal M}_{\infty} \rightarrow \grO^3_0 SU(2)$ 
be the
direct limit of the maps $\vartheta_k$ in diagram 
{\rm(1.1)}. Then
$\vartheta$ is a homotopy equivalence.   
\ethm

Using different methods Gravesen [Gra] independently proved 
that $\vartheta$ is a homology equivalence.  However, these 
are stable results and do not directly address the 
structure 
of $H_t({\cal M}_k)$ or $\pi_t({\cal M}_k)$ for any fixed 
$t$ 
or $k$.  In another direction, [BoMa] showed that 
$(\vartheta_k)_t$ is a surjection in $\bbz/p$ homology for 
$t\leq k$ but the techniques used there do not extend to 
analyze the kernel.

In [BHMM1] we prove the following theorem, which settles 
conjectures
1, 2, and 3 of Atiyah and Jones in the affirmative.

\thm{Theorem A}
For all $k>0$, the natural inclusion 
$\vartheta_k:{\cal M}_k \rightarrow \grO^3_k SU(2)$ 
is a homotopy equivalence through dimension 
$q=q(k)= [k/2]-2$.  
\ethm

By a homotopy equivalence through dimension $q(k)$ we mean 
that
$\vartheta_k$ induces an isomorphism $\pi_t({\cal M}_k) 
~\cong~
\pi_t(\grO^3_k SU(2))= \pi_{t+3}(SU(2))$ for $t\leq q(k)$. 
 Thus,
for example, the low-dimensional homotopy groups of ${\cal 
M}_k$
are finite.   The statement in homotopy implies the weaker
statement that $\vartheta_k$ induces an isomorphism 
$H_t({\cal M}_k;A) ~\cong~ H_t(\grO^3_k SU(2);A)$
for $t\leq q(k)$ and any coefficient group $A$.  Here, the 
homology
groups on the right are completely known. 

Actually, Theorem A is a corollary of the main topological 
result
in [BHMM1].

\thm{Theorem B}
For all $k>0$, the 
inclusion $\gri_k:{\cal M}_k \rightarrow {\cal M}_{k+1}$ 
is a homotopy
equivalence through dimension 
$q=q(k)= [k/2]-2$.  
\ethm

These topological results are actually formal consequences 
of our
detailed geometric analysis of ${\cal M}_k$. This is 
explained in
the following section.  Lastly, it is likely that the 
explicit
bound $q(k)$ can be improved as explained at the end of 
the next
section.

\heading 2. Outline of proofs\endheading

A key tool in our approach is a theorem of Donaldson [D], 
that
allows us to work in a purely holomorphic context and 
study a
certain moduli space of holomorphic bundles.  To prove 
Theorems A
and B we study this equivalent holomorphic moduli space by
following the program used in [MM1]  and [MM2] to study the
topology of holomorphic maps from the Riemann sphere into 
flag
manifolds.  This should not be surprising since these 
holomorphic
mapping spaces correspond, by results of [Hu2] and [HM], 
to moduli
spaces for $SU(n)$-monopoles. The main idea is to 
construct a
stratification of the entire moduli space by smooth 
manifolds, each
with an orientable normal bundle, and each contained in 
the closure
of higher dimensional strata. These strata are then 
organized into
a filtration, and the resulting spectral sequence is 
analyzed to
prove homology versions of Theorems  A and B.  More 
precisely, we
prove

\thm{Theorem C}
For each $k$ and all
coefficient rings $A$, there are homology Leray spectral 
sequences
$E^r({\cal M}_k;A)$ converging to filtrations of 
$H_*({\cal M}_k;A)$ with identifiable $E^1$ terms \RM(see\ 
\cite{BHMM1}
for the precise statement\RM).  Furthermore, the inclusion 
$\gri_k: {\cal M}_k \to {\cal M}_{k+1}$ induces a map of
spectral sequences $\gri_k(r)\: E^r({\cal M}_k;A) \to 
E^r({\cal M}_{k+1};A).$ In particular, differentials in 
these
spectral sequences are natural with respect to $\gri_k(r)$.
\ethm

Next, we show that the map $\gri_k$ induces an
isomorphism of the $E^2$ terms 
in Theorem C through dimension $q(k)+2$.   Since 
differentials are
natural, it follows that $\gri_k$ induces a homology 
equivalence
through dimension $q(k)+1$.  This gives the homology 
version of
Theorem B, namely,
\thm{Theorem D}
For all $k>0$ and all
coefficients $A$, the inclusion
$\gri_k\: {\cal M}_k \to {\cal M}_{k+1}$ 
induces an isomorphism in homology 
$(\gri_k)_t: H_t({\cal M}_k;A) ~\cong~ H_t({\cal M}_{k+
1};A)$ 
for $t\leq q=q(k)+1.$ 
\ethm

Combined with the homology Taubes-Gravesen 
stability theorem, Theorem D 
implies the homology version of Theorem A.  However, as 
${\cal M}_k$ and $\grO^3 S^3$ are not simply connected,
Theorems A and B do not trivially follow from their 
homology 
analogs.  Thus, we analyze the induced map on the
universal covers $(\tilde{\vartheta}_k)_t: 
H_*(\widetilde{\cal M}_k;\bbz) \to 
H_t(\widetilde{\grO}^3 S^3;\bbz)$
and show this is an isomorphism for $t\leq q(k)$. 
This now implies Theorems A and B.

Thus, to complete the outline of the proof we need to 
explain how
we prove Theorem C.  It is nontrivial to obtain a 
stratification
for ${\cal M}_k$ that is suited to this type of topological
analysis. In fact, much of the work in [BHMM1] is devoted to
understanding the geometry of ${\cal M}_k$ in sufficient 
detail so
that we can construct such a stratification.  As mentioned 
above
the starting point for this analysis is a result of 
Donaldson [D] 
saying that the framed $SU(2)$ instanton moduli space is 
equivalent 
to the moduli space of rank two holomorphic vector bundles 
on $\bbp^2$ 
with $c_1=0$ that are trivial when restricted to a fixed 
line $\ell$ 
with a fixed trivialization there.
Moreover, as was noticed in [A] and [Hu1], it is 
convenient to perform a
birational transformation on the line $\ell$ to obtain a 
surface ruled by
lines, for example, $\bbp^1\times \bbp^1$.  Let ${\cal 
M}_k(M,S)$ denote the
moduli space of isomorphism classes of rank 2 holomorphic 
bundles on $M$ with
$c_1=0$ and $c_2=k$ that have a fixed trivialization on 
$S$.  Then there are
diffeomorphisms
$${\cal M}_k\simeq {\cal M}_k(\bbp^2,\ell)\simeq 
{\cal M}_k(\bbp^1\times\bbp^1,{\ell}_1\vee 
{\ell}_2).$$

Thus, the geometry of ${\cal M}_k$ is described by 
isomorphism
classes of certain framed semistable holomorphic rank 2 
vector
bundles $E$ over $\bbp^1\times\bbp^1$.  It is known that 
such $E$
are trivial on almost all lines $\bbp^1\times \{x\}$ and 
that the structure
of $E$ is determined by its behaviour on neighbourhoods of 
a finite number
of {\it jumping lines} $\bbp^1\times\{x_i\}$, lines over 
which the holomorphic
structure of E is that of a sum of line bundles ${\cal 
O}(h)\oplus{\cal O}(-h)$
with $h>0$.  In particular, if $\{\infty\}$ in the second 
$\bbp^1$ 
corresponds to the
line $\ell_1$, then the possible jumping lines are 
parameterized by the
complex plane $\bbc\simeq \bbp^1-\{\infty\}$.  At each 
jumping line there is an
associated integer $m$ called the multiplicity.  Let $Q_m$ 
denote the framed
isomorphism classes (on a neighbourhood of a jumping line) 
of framed jumps of
multiplicity $m$.   To each instanton $\gro$ is thus 
associated an isomorphism
class $[E_{\gro}]$ of bundles $E_{\gro}$ and, in turn, to 
each $[E_{\gro}]$ is
associated a point in a labelled configuration space. 
Elements of 
this space are given by  configurations of points
$z_i$ in the complex plane $\bbc$ determined by a finite 
number of nontrivial
jumping lines each labelled by elements $l_i\in Q_{m_i}$.  
Here the total
charge $k$ is the sum of the multiplicities.  Thus, as 
shown in [Hu1], there is
a natural holomorphic projection $\Pi:{\cal M}_k 
\to SP^k(\bbc)$ that
associates to any equivalence class of framed bundles 
$[E]$ its divisor of
jumping lines $\sum_{i=1}^r m_iz_i$ in $SP^k(\bbc)$.  The 
fibre at a point
$z\in SP^k(\bbc)$ with multiplicity $(m_1,\dots,m_r)$ is 
the product
$Q_{m_1}\times \cdots \times Q_{m_r}$ where $Q_m$ denotes 
the space of
equivalence classes of framed jumps of multiplicity $m$. 
This is the ``pole and
principal part" picture of Segal-Gravesen [Gra].

However, as each $Q_m$ is not necessarily a smooth 
manifold, this decomposition
of ${\cal M}_k$ is yet not sufficiently fine to apply the
techniques of [MM1, MM2].  Therefore, a refined local 
analysis of
$Q_m$ is required.  To each multiplicity $m_i$ we associate 
finite sequences $G=(h_0,h_1,\ldots,h_{l-1})$ of decreasing 
integers called ``heights" such that 
$\sum^{l-1}_{j=0}h_j=m_i$.  
We refer to these sequences   $G$ as ``graphs"
and show that to each framed jumping line of multiplicity 
$m$ we
can associate a unique graph $G$. This is done by setting, 
for each $k$,
$$
h_k = \dim\left\{\matrix
\hbox{\rm sections\ of\ }E\otimes {\cal O}(-1) 
\hbox{\rm\  on\ the\ jumping\ line\ which\  extend }\\
\hbox{\rm\ to\ the\ } k\hbox{\rm th\ formal\ 
neighbourhood\ of\ the\ line\ in\ }
 \bbp^1\times\bbp^1\endmatrix
\right\}. 
$$ 
We thus  obtain a decomposition 
$Q_m = \bigsqcup_G FJ_G$ where each $FJ_G$ is the 
collection of
framed jumps of a fixed graph type and the union is taken 
over all
graphs $G$ of multiplicity $m$.  Most importantly, we prove
\thm{Theorem E}
Let $G$ be a graph with heights
$h=h_0\geq h_1\geq\cdots\geq h_{l-1}>h_l=0$, length $l$, and
multiplicity $m=\sum h_i$.  The space of framed jumps 
$FJ_G$ with
graph $G$ is a smooth complex manifold of complex 
dimension $2m+l$.
\ethm

Let $S_{G_1,\dots,G_r}$ denote the set of equivalence 
classes
of framed holomorphic bundles with $r$ jumping lines
$(L_1,\dots,L_r)$ with graphs ${\cal G}(r)=(G_1,\dots,G_r)$,
respectively. Thus, we have the disjoint union decomposition
${\cal M}_k = \bigsqcup S_{{\cal G}}$. Furthermore, the 
map $\Pi$
now restricts to a map $\Pi_{{\cal G}(r)}\:S_{{\cal G}(r)} 
\to\bbd\bbp^r(\bbc)$ with fibre 
$FJ_{G_1}\times\cdots\times FJ_{G_r}$.
Here $\bbd\bbp^r(\bbc)$ is the deleted product; that is, 
the space
of $r$ distinct unordered points in $\bbc$, which is a 
smooth
complex manifold of complex dimension $r$.  From Theorem E 
follows

\thm{Theorem F}
$S_{{\cal G}(r)}$
is a smooth complex variety of complex dimension $2k+l+r$. 
\ethm

We then prove that this decomposition of ${\cal M}_k$ by the
$S_{{\cal G}}$ submanifolds is a stratification of the type
required.  In particular, it is necessary to know the 
codimension
of each stratum (which is given in the last theorem), that 
there
is a natural lexicographical order on the graph types, 
and, with
respect to that order, what the intersection of the normal 
bundle
of a fixed stratum is with the other strata.  In addition to
checking these facts we verify that our constructions are 
all
compatible with the stabilization maps $\gri_k$ in diagram
(1.1).  After combining these results Theorem C follows as 
in
[MM1, MM2].

Finally, we note that the stability argument given here 
does not
require explicit calculation of differentials in the 
spectral
sequence given in Theorem C but uses only information at 
the $E^1$
level and naturality.  In a sequel [BHMM2] we will examine 
the
structure of the $E^1$ term and of the differentials much 
more
carefully.  This will permit us to both sharpen the 
explicit value
of $q(k)$ and, more importantly, obtain much more detailed
information about $H_*({\cal M}_k)$ above the range of 
stability.



\Refs\ra\key{BHMM2}

\ref\key A 
\by M. F. Atiyah \paper Instantons in two and four
dimensions  
\jour Comm. Math. Phys. \vol 93 
\yr 1984  \pages 437--451
\endref

\ref\key AJ 
\by M. F. Atiyah and J. D. Jones \paper Topological
aspects of Yang-Mills theory  
\jour Comm. Math. Phys. \vol 61 \yr 1978 
\pages 97--118
\endref

\ref\key BHMM1 
\by C. P. Boyer, J. C. Hurtubise, B. M. Mann, and R. J.
Milgram \paper The topology of instanton moduli spaces. 
{\rm I:} The
Atiyah-Jones conjecture  
\jour Ann. of Math. (2), \toappear 
\endref

\ref\key BHMM2 \bysame  
\paper The topology of instanton moduli spaces. {\rm II:} 
The
Toeplitz varieties  
\paperinfo in preparation
\endref

\ref\key BoMa 
\by C. P. Boyer and B. M. Mann \paper Homology
operations on instantons  
\jour J. Differential Geom. \vol 28 
\yr 1988  \pages 423--465
\endref

\ref\key D 
\by S. K. Donaldson \paper Instantons and geometric
invariant theory  
\jour Comm. Math. Phys. \vol 93 
\yr 1984 \pages 453--461
\endref

\ref\key Gra 
\by J. Gravesen \paper On the topology of spaces
of holomorphic maps  
\jour Acta Math. \vol 162 
\yr 1989 \pages 247--286
\endref

\ref\key Hu1 
\by J. C. Hurtubise \paper Instantons and jumping
lines  
\jour Comm. Math. Phys. \vol 105 
\yr 1986  \pages 107--122
\endref

\ref\key Hu2 \bysame \paper The classification of
monopoles for the classical groups  
\jour Comm. Math. Phys. \vol 120 
\yr 1989  \pages 613--641
\endref

\ref\key HM 
\by J. C. Hurtubise and M. K. Murray \paper On the
construction of monopoles for the classical groups  
\jour Comm. Math. Phys. \vol 122
\yr 1989  \pages 35--89
\endref

\ref\key MM1 
\by B. M. Mann and R. J. Milgram \paper Some spaces
of holomorphic maps to complex Grassmann manifolds  
\jour J. Differential Geom. \vol 33 
\yr 1991 \pages 301--324
\endref

\ref\key MM2 \bysame  \paper On the moduli
space of SU(n) monopoles and holomorphic maps to flag 
manifolds  
\paperinfo preprint, UNM and Stanford University, 1991
\endref

\ref\key T1 
\by C. H. Taubes \paper Path-connected Yang-Mills
moduli spaces  
\jour J. Differential Geom. \vol 19 
\yr 1984 \pages 337--392
\endref

\ref\key T2 \bysame \paper The stable topology of
self-dual moduli spaces  
\jour J. Differential  Geom. \vol 29
\yr 1989 \pages 163--230
\endref
\endRefs
\enddocument